\documentclass{amsart}

\usepackage[all]{xy}

\usepackage{pictexwd, epsfig, amsthm, amssymb, amsmath, graphicx, psfrag,amsxtra, latexsym}
\usepackage{enumerate}

\newcommand{\mD}{\mathbb D}
\newcommand{\mZ}{\mathbb Z}
\newcommand{\mR}{\mathbb R}

\newcommand{\mT}{\mathcal T}

\def\polhk#1{\setbox0=\hbox{#1}{\ooalign{\hidewidth
    \lower1.5ex\hbox{`}\hidewidth\crcr\unhbox0}}}

\theoremstyle{plain}
\newtheorem{Thm}{Theorem}
\newtheorem{Cor}[Thm]{Corollary}
\newtheorem{Prop}[Thm]{Proposition}

\newtheorem{lemma}[Thm]{Lemma}

\theoremstyle{definition}
\newtheorem{definition}[Thm]{Definition}

\theoremstyle{remark}

\begin{document}

\title[CAT(0)-manifolds with no Riemannian smoothings]
      {4-dimensional locally CAT(0)-manifolds with no Riemannian smoothings}

\author{M. Davis}
\address{Department of Mathematics\\
         Ohio State University\\
         Columbus, OH  43210}
\email[Davis]{mdavis@math.ohio-state.edu}

\author{T. Januszkiewicz}
\address{Department of Mathematics\\
         Ohio State University\\
         Columbus, OH  43210}
\email[Januszkiewicz]{tjan@math.ohio-state.edu}

\author{J.-F. Lafont}
\address{Department of Mathematics\\
         Ohio State University\\
         Columbus, OH  43210}
\email[Lafont]{jlafont@math.ohio-state.edu}

\begin{abstract}
We construct examples of 4-dimensional manifolds $M$ supporting a locally
CAT(0)-metric, whose universal covers $\tilde M$ satisfy Hruska's isolated flats condition, and 
contain 2-dimensional flats $F$ with the property that $\partial ^\infty F \cong S^1 \hookrightarrow 
S^3 \cong \partial ^\infty \tilde M$ are nontrivial knots. As a consequence, we 
obtain that the group $\pi_1(M)$ cannot be isomorphic to the fundamental
group of any Riemannian manifold of nonpositive sectional curvature. In particular, if $K$
is any locally CAT(0)-manifold, then 
$M\times K$ is a locally CAT(0)-manifold which does not support any Riemannian metric of nonpositive
sectional curvature.

\end{abstract}

\maketitle

\section{Introduction}

Riemannian manifolds of nonpositive sectional curvature are a class of manifolds
featuring a rich interplay between their geometry, their topology, and their dynamics. In the 
broader setting of geodesic metric spaces, we have the notion of a locally CAT(0)-metric.
These provide a metric space analogue of nonpositively curved Riemannian manifolds,
and many classic results concerning Riemannian manifolds of nonpositive sectional curvature have
now been shown to hold more generally for locally CAT(0)-spaces. We are interested in 
understanding the difference, within the class of closed manifolds, between (1) supporting
a Riemannian metric of nonpositive sectional curvature, and (2) supporting a locally 
CAT(0) metric. A closed topological manifolds equipped with a locally CAT(0)-metric will be called
a {\it locally CAT(0)-manifold}.

In low dimensions, there is no difference between these two classes. In two dimensions, 
this follows easily from the classification of surfaces, while in three dimensions, this follows
from Thurston's geometrization theorem (recently established by Perelman). In contrast, Davis
and Januszkiewicz \cite{DJ} have constructed examples, in all dimensions $\geq 5$, of locally CAT(0)-manifolds 
which do {\it not} support any Riemannian metric of nonpositive
sectional curvature. In this paper, we deal with the remaining open case.

\vskip 10pt

\noindent {\bf Main Theorem:} There exists a 4-dimensional closed manifold $M$ with
the following four properties:
\begin{enumerate}
\item $M$ supports a locally CAT(0)-metric, 
\item $M$ is smoothable, and $\tilde M$ is diffeomorphic to $\mR ^4$,
\item $\pi_1(M)$ is {\bf not} isomorphic to the fundamental group of 
any Riemannian manifold of nonpositive sectional curvature.
\item if $K$ is any locally CAT(0)-manifold, then $M\times K$ is a locally CAT(0)-manifold which does not support any 
Riemannian metric of nonpositive sectional curvature.
\end{enumerate}

\vskip 15pt

Let us briefly outline the idea behind the proof of our main result. First of all, we introduce the
notion of a triangulation of $S^3$ to have {\it isolated squares}. Any such triangulation has a well-defined
{\it type}, which is the isotopy class of an associated link in $S^3$. In Section 3, we provide a proof that any given
link in $S^3$ can be realized as the type of a suitable 
flag triangulation of $S^3$ with isolated squares. In Section 4, we start with a flag triangulation $L$
of $S^3$ with isolated squares, whose type is a nontrivial knot, and use it to construct the desired $4$-manifold.
This is done by considering
the right angled Coxeter group $\Gamma_{L}$ associated to the triangulation $L$, and defining $M$ to be the 
quotient of the corresponding Davis complex by a torsion free finite index subgroup $\Gamma \leq \Gamma_L$.
Standard properties of the triangulation $L$ ensure that $M$ is smoothable, and that the Davis complex is CAT(0) and
diffeomorphic to $\mR^4$. The isolated squares condition on the flag triangulation $L$ ensures the Davis complex satisfies
Hruska's {\it isolated flats} condition. The fact that the type of $L$ is a nontrivial knot ensures that the Davis complex
contains a periodic $2$-dimensional flat $F$ which is {\it knotted at infinity}. But now if $M$ supported a Riemannian metric $g$ of 
nonpositive sectional curvature, the flat torus theorem ensures that one could find a corresponding flat $F^\prime$ (in the $g$-metric) 
which is $\Gamma$-equivariantly homotopic to $F$, and the isolated flats condition then forces $F^\prime$ to also be knotted at infinity. 
However, in the Riemannian setting, it is easy to see that a codimension two flat must be unknotted at infinity, yielding a contradiction.

\vskip 10pt

\centerline{\bf Acknowledgments}

\vskip 5pt

The first two authors were partially supported by the NSF, under grant DMS-50706259.
The last author was partially supported by the NSF, under grant DMS-0906483, and by an Alfred P.
Sloan research fellowship.

\vskip 10pt

\section{Previously known obstructions.}

Our Main Theorem provides a new obstruction to the problem of finding a {\it Riemannian smoothing} 
on a manifold $M$ supporting a locally CAT(0)-metric. More precisely, we say that such a manifold supports
a Riemannian smoothing provided one can find a smooth Riemannian manifold $(N,g)$, with $g$ a 
Riemannian metric of nonpositive sectional curvature, and a homeomorphism $f: N \rightarrow M$.
In this section, we briefly summarize the known obstructions to Riemannian smoothing.

\subsection{Example: no smooth structure.} Given a Riemannian smoothing
$f: N \rightarrow M$ of a locally CAT(0)-manifold $M$, one can forget the Riemannian structure and simply
view $N$ as a smooth manifold. This immediately tells us that, if $M$ has a Riemannian smoothing, then it
must be homeomorphic to a smooth manifold, i.e. the topological manifold $M$ must be {\it smoothable}. The 
first examples of aspherical topological manifolds not homotopy equivalent to smooth manifolds were 
constructed (in all dimensions $\geq 13$) by Davis and Hausmann \cite{DH} by using the reflection group trick. 
Non-smoothable aspherical PL-manifolds were constructed (in all dimensions $\geq 8$)
in the same paper. For the sake of completeness, we now sketch out a (slightly different) construction of a
closed 8-dimensional locally CAT(-1)-manifold $M^8$ which is 
not homotopy equivalent to any smooth 8-manifold. 

Recall that Milnor constructed \cite{Mi} an 8-dimensional PL-manifold $N^8$ which
is not homotopy equivalent to any smooth 8-manifold. Milnor's example had the property that the second rational 
Pontrjagin class $p_2(N^8)$ was {\it not} an integral class, and hence cannot be homeomorphic to a smooth manifold. 
Let us take $N^8$ equipped with a PL-triangulation. Charney and Davis \cite{CD} developed a {\it strict hyperbolization} 
process, which inputs a triangulated manifold $M$ and outputs a piecewise hyperbolic manifold $h(M)$ 
equipped with a locally CAT(-1)-metric. Furthermore, they showed that the hyperbolization process 
preserves rational Pontrjagin classes. In particular, applying their strict hyperbolization process to $N^8$,
we obtain a locally CAT(-1)-manifold $h(N^8)$, having the property that $p_2(h(N^8))$ fails to be integral, and
hence forcing $h(N^8)$ to be non-smoothable. Finally, we note that the Borel Conjecture is known to hold for
this class of aspherical manifolds (see \cite{BL}), so if $h(N^8)$ was homotopy equivalent to some smooth manifold,
it would in fact be homeomorphic to the smooth manifold (contradicting non-smoothability). Similar examples
can be constructed in all dimensions of the form $n=4k$, with $k\geq 2$ (see also the discussion in
\cite[Section 5]{BLW}).

\subsection{Example: no PL structure.}\label{ss:PL}  
In a similar vein, it is also possible to construct (topological) locally CAT(0)-manifolds that do not even support any PL-structures. 
We recall such an example from \cite[Section 5a]{DJ}.  We let $M^4(E_8)$ denote the $E_8$ homology manifold. Recall that
this space is constructed by first plumbing together eight copies of the tangent disk bundle to $S^2$, according to the pattern given
by the $E_8$ Dynkin diagram. This results in a smooth $4$-manifold with boundary $N^4$, whose boundary $\partial N^4$ is 
homeomorphic to
Poincar\'e's homology $3$-sphere. Coning off the boundary gives the space $M^4(E_8)$, a simply connected homology manifold 
of signature $8$ with one singular point. Taking a triangulation of $N^4$, one can extend it (by coning on the boundary) 
to a triangulation of $M^4(E_8)$, which we can then hyperbolize to obtain a space $H^4$.

The space $H^4$ is now a homology $4$-manifold of signature $8$ with one singular point, and comes equipped with a locally 
CAT(0)-metric.  It follows from Edward's Double Suspension Theorem that $H^4\times T^k$ is a topological $(4+k)$-manifold
(where $T^k$ denotes the $k$-torus and $k\geq 1$). The manifolds $H^4 \times T^k$ come equipped with a 
(product) locally CAT(0)-metric, but it follows from the arguments in \cite[Section 5a]{DJ} that they do not admit 
a PL structure. Thus, in each dimension $\geq 5$ there is a locally CAT(0)-manifold with no PL structure.

\subsection{Example: universal cover distinct from $\mR ^n$.} For a third family of examples, we 
recall that the classic Cartan-Hadamard theorem asserts that the universal cover of a Riemannian
manifold of nonpositive sectional curvature must be diffeomorphic to $\mathbb R^n$. In particular,
a CAT(0)-manifold $M$ with the property that $\tilde M$ is {\it not} diffeomorphic to $\mathbb R^n$
can not support a Riemannian smoothing. Davis and Januszkiewicz constructed (see \cite[Thm. 5b.1]{DJ}) 
examples of locally CAT(0)-manifolds $M^n$ (for $n\geq 5$), with the property 
that their universal covers $\tilde M^n$ are {\it not} simply connected at infinity (and hence, not 
homeomorphic to $\mR^n$). Further examples of this type are described in \cite{adg}.

\subsection{Example: boundary at infinity distinct from $S^{n-1}$.} In the previous three families of examples,
{\it topological} properties (smoothability, PL-smoothings, topology of universal cover) were used to obstruct the existence
of a Riemannian metric of nonpositive sectional curvature. The next family of examples have obstructions
that arise from the {\it large scale geometry} of the universal covers. Associated to a CAT(0)-space $X$, we have
a topological space called the {\it boundary at infinity} $\partial ^\infty X$.  If $X$ is Gromov hyperbolic, then the 
homeomorphism type of $\partial^\infty X$ is a quasi-isometry invariant of $X$.  In particular, if $X$ is the 
universal cover of a locally CAT(-1)-space $Y$, then $\partial^\infty X$ depends only on $\pi_1(Y)$.  When $X$ is the 
universal cover of an $n$-dimensional closed Riemannian manifold of nonpositive sectional curvature, the
corresponding $\partial ^\infty X$ is homeomorphic to the standard sphere $S^{n-1}$.  

Now consider the locally CAT(-1) $5$-manifold $M^5$ obtained by applying a strict hyperbolization procedure 
(from \cite{CD}) to the double suspension of a triangulation of Poincar\'e's homology $3$-sphere.
Denote by $X^5$ its universal cover, and observe that, although $\partial^\infty X^5$ has the homotopy type of $S^4$, 
it is proved in \cite[Section 5c]{DJ} that $\partial ^\infty X^5$ is {\it not} locally simply connected. So 
$\partial ^\infty X^5$ cannot be homeomorphic to $S^4$ (in fact, is not even an ANR).  Thus, $M^5$ is not 
homotopy equivalent to a Riemannian $5$-manifold of strictly negative sectional curvature.  The same argument 
applies to a strict hyperbolization of the manifold $M^4(E_8)\times S^1$ discussed in Section \ref{ss:PL}.  There are 
similar examples in higher dimensions $n>5$ obtained by strictly hyperbolizing double suspensions of homology 
$(n-2)$-spheres.  Thus, in each dimension $n \geq 5$ there are closed locally CAT(-1) manifolds $M^n$ with 
universal cover homeomorphic to $\mR^n$ but which are not homotopy equivalent to any Riemannian $n$-manifold 
of strictly negative sectional curvature.

\subsection{Example: stability under products.} Finally, we point out one last method for producing
manifolds which do not have Riemannian smoothings:

\begin{Prop}
Let $M^n$ be a locally CAT(0)-manifold which does not support any Riemannian smoothing, and assume that $n\geq 5$.
Then for $K$ an arbitrary locally CAT(0)-manifold, the product $M\times K$ is a locally CAT(0)-manifold which does not support
any Riemannian smoothing.
\end{Prop}

\begin{proof}
To see this, we first note that the product of the locally CAT(0)-metrics on $M$ and $K$ provide a locally CAT(0)-metric 
on $M\times K$. Now assume that $M\times K$ supported a Riemannian smoothing $f: N\rightarrow M\times K$, 
and let $g$ be the associated Riemannian metric of nonpositive sectional curvature on $N$. 
Since $\pi_1(N) \cong \pi_1(M)
\times \pi_1(K)$, the classical
splitting theorems (see Gromoll and Wolf \cite{GW}, Lawson and Yau \cite{LY}, and Schroeder \cite{Sc})
imply that we have a corresponding {\it geometric} splitting $(N, g) \cong (M^\prime, g_1) \times 
(K^\prime, g_2)$, having the property that:
\begin{itemize}
\item each factor can be identified with a totally geodesic submanifold of 
$(N, g)$, 
\item the factors satisfy $\pi_1(M)\cong \pi_1(M^\prime)$, and $\pi_1(K) \cong \pi_1(K^\prime)$.
\end{itemize}
So we see that $M^\prime$ is a Riemannian manifold of nonpositive sectional curvature, of
dimension $\geq 5$, and satisfying $\pi_1(M) \cong \pi_1(M^\prime)$. Since the Borel conjecture is known
to hold for this class of manifolds (see Farrell and Jones \cite{FJ}), there exists a homeomorphism 
$M^\prime \rightarrow M$ realizing the isomorphism of fundamental groups. This provides a 
Riemannian smoothing of $M$, giving us the desired contradiction.
\end{proof}

We remark that property (4) in our Main Theorem can be deduced from a virtually identical argument: instead
of appealing to the Borel Conjecture to obtain a contradiction, we resort instead to property (3) in our Main Theorem. 

\section{Special triangulations of $S^3$.} 

Recall that a simplicial complex is {\it flag} provided it is determined
by its 1-skeleton, i.e. every $k$-tuple of pairwise incident vertices spans a $(k-1)$-simplex $\sigma ^{k-1}$
(for $k\geq 3$). A subcomplex $\Sigma^\prime$ of a simplicial complex $\Sigma$ is {\it full} provided every
simplex $\sigma \subset \Sigma$ whose vertices lie in $\Sigma ^\prime$ satisfies $\sigma \subset \Sigma^\prime$.
We will say a cyclically ordered 4-tuple of vertices $(v_1, v_2, v_3, v_4)$ in a simplicial complex forms a {\it square} 
provided each consecutive pair of vertices determines an edge in the complex, while the pairs $(v_1,v_3)$ 
and $(v_2,v_4)$ do {\it not} determine an edge.

\begin{figure}
\label{graph}
\begin{center}
\includegraphics[width=2in, angle=0]{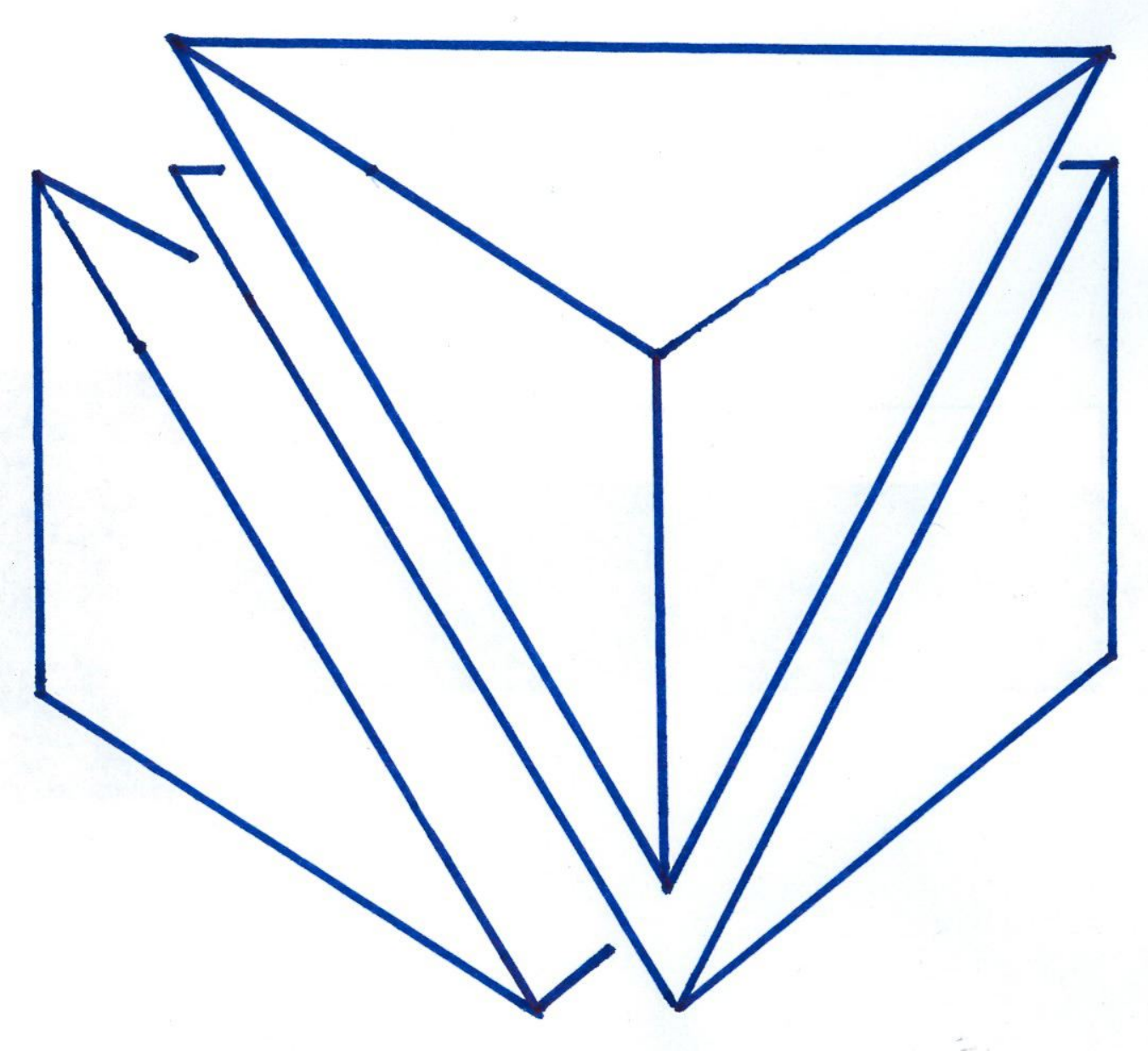}
\caption{Basic triangulation of a triangular prism.}
\end{center}
\end{figure}

\begin{definition}
A flag triangulation of $S^3$ is said to have {\em isolated squares} provided no two squares in the triangulation
intersect (i.e. each vertex lies in at most one square). For such a triangulation, the collection of squares form
a link in $S^3$. We call the isotopy class of this link the {\em type} of the triangulation.
\end{definition}

In this section, we establish:

\begin{Thm}
Let $k\subset S^3$ be any prescribed link in the 3-sphere. Then there exists a flag triangulation of $S^3$,
with isolated squares, and with type the given link $k$.
\end{Thm}

We establish this result in several steps, gradually building up the triangulation to have the properties we desire. 

\vskip 10pt

\noindent{\bf Step 1: Triangulating the solid torus.}

\vskip 5pt

As a first step, we describe a triangulation on a solid torus $\mD^2 \times S^1$. Recall that there is 
a canonical decomposition of the 3-dimensional cube $[0,1]^3 \subset \mR^3$ into six tetrahedra. This triangulation is 
determined by the inequalities $0\leq x_{\sigma(1)} \leq x_{\sigma(2)}\leq x_{\sigma (3)}\leq 1$, where $\sigma$ 
ranges over the six possible permutations of the index set $\{1,2,3\}$. Now if we restrict to the region where 
$x_1\leq x_2$, we obtain a  triangulation of the triangular prism $\Delta ^2 \times [0,1]$ into
exactly three tetrahedra. Let us denote by $F, G$ the two square faces of the triangular prism defined via the
hyperplanes $x_1=0$ and $x_1=x_2$ respectively. The triangulation of the prism cuts each of these squares
into two triangles, along the diagonal originating at the origin. We call the {\it bottom} of the prism the triangle
corresponding to the intersection with the hyperplane $x_3=0$, and call the {\it top}
of the prism the triangle arising from the intersection with the hyperplane $x_3 = 1$. Figure 1 contains an illustration
of this decomposition of the triangular prism (drawn to respect the orientation of the ``bottom'' and ``top''). In the picture,
the two square sides facing us are $F$ and $G$ respectively. 

We can now take three copies of the triangular prism, and cyclically identify each $F_i$ to the 
corresponding $G_{i+1}$. This gives a new triangulation of a triangular prism (with nine tetrahedra), 
with an inherited notion of ``top'' and ``bottom''. This new triangulation has the following key properties:
\begin{itemize}
\item there exists a unique edge $e$ of the triangulation joining the center of the bottom triangle to the center of the top 
triangle,
\item the center of the bottom triangle is adjacent to {\it every} vertex in the triangulation, and
\item aside from the center of the bottom triangle, the center of the top triangle is adjacent to {\it no other}
vertices in the bottom of the prism.
\end{itemize}
We will call a copy of this canonical triangulation of the triangular prism a {\it block}. Fixing an
identification of $\mD^2$ with the base of the triangular prism, we can think of a block as a 
triangulation of $\mD^2\times [0,1]$. 

To obtain the desired triangulation of the solid torus $\mD^2\times S^1$, we ``stack'' four blocks together. 
More precisely, we take four blocks and cyclically identify the top of each block with the bottom of the next block. 
This gives us a triangulation of the solid torus $\mD^2 \times S^1$ into thirty-six tetrahedra. We
say blocks are {\it adjacent} or {\it opposite}, according to whether they share a vertex or not. Corresponding
to the above properties for the individual blocks, this triangulation of the solid torus
satisfies:
\begin{itemize}
\item the triangulation contains a canonical, unique square having the property that it is entirely contained
within the {\it interior} of $\mD^2\times S^1$; the four vertices of this square will be called {\it interior vertices}.
\item all the remaining vertices of the triangulation lie on the boundary of $\mD^2\times S^1$, and will be
called {\it boundary vertices}.
\item every tetrahedron in the triangulation contains at least one interior vertex.
\item every interior vertex has the property that, if one looks at all adjacent boundary vertices, these vertices 
are all contained in single block (the unique block whose bottom contains the given interior vertex).
\end{itemize}
We call the unique square in the interior of this triangulation of $\mD^2\times S^1$ the {\it core} of the solid torus.
Observe that, out of the thirty-six tetrahedra occuring in the triangulation, exactly twenty-four of them arise
as the join of a triangle in $\partial \mD^2\times S^1$ with an interior vertex, while the remaining twelve occur as the
join of an edge in $\partial \mD^2\times S^1$ with an edge in the core.

\vskip 10pt

\noindent {\bf Step 2: Getting squares realizing the link $k$.}

\vskip 5pt

Next, let us take the desired link $k$, and take pairwise disjoint regular closed neighborhoods $\hat N_i$ of the individual 
components of the link. Each of these neighborhoods is homeomorphic to a solid torus, and we denote by 
$N_i \subset \hat N_i$ the slightly smaller solid torus of radius half as large. We proceed to construct a 
triangulation of $S^3$ as follows: first, within each of the tori $N_i$, we use the triangulation described
in Step 1, identifying the components of the link with the cores of the various triangulated solid tori. 
Secondly, removing the interiors of all of the $\hat N_i$, we obtain a compact 3-manifold $M$ with 
boundary $\partial M = \coprod \partial \hat N_i$. Since 3-manifolds are triangularizable, we now choose an arbitrary 
triangulation of this 3-manifold $M$, obtaining a triangulation of $M \cup \coprod N_i \subset S^3$.
The closure of the complementary region is a disjoint union of the sets $\hat N_i \setminus N_i$, each of 
which is topologically a fattened torus $S^1 \times S^1 \times [0,1]$. Furthermore, we are given triangulations 
$\mathcal T_0, \mathcal T_1$ of the two
boundaries $S^1\times S^1 \times \{0\}$, $S^1\times S^1\times \{1\}$ (coming from the triangulations 
of $\partial N_i$ and $\partial M$ respectively). But any two triangulations of the 2-torus $S^1\times S^1$ have
subdivisions which are simplicially isomorphic. Letting $\mathcal T ^\prime$ denote such a triangulation, we
assign this triangulation on the level set $S^1\times S^1 \times \{1/2\}$. 

Finally, we extend the triangulation into the two regions $S^1 \times S^1 \times [0,1/2]$ and 
$S^1 \times S^1 \times [1/2,1]$ using the following procedure. On each of these two regions, we have
a triangulation $\mT_i$ on one of the boundary components, and a subdivision $\mT ^\prime$
of the triangulation on the other boundary component. We proceed to inductively subdivide each of the
regions $\sigma \times I$, where $\sigma$ ranges over the simplices of the triangulation $\mT_i$.
First of all, we add in edges $\sigma^0\times I$ for each vertex in the triangulation $\mT_i$. Now assuming
that we have already triangulated the product $\mT_i^{(k-1)}\times I$ of the $(k-1)$-skeleton of $\mT_i$ with 
the interval, let us extend the triangulation to $\mT_i^k\times I$. Given a $k$-simplex 
$\sigma^k$, we have that the region $\sigma^k \times I$ is topologically a closed $(k+1)$-dimensional
ball, with boundary that can be identified with $(\sigma^k \times \{0\}) \coprod (\sigma ^k\times \{1\}) \coprod
(\partial \sigma ^k \times I)$. Furthermore, the bottom level consists of a simplex (the original $\sigma ^k \in \mT_i$), 
the top level consists of a subdivision of the simplex (the subdivision of $\sigma ^k$ inside $\mT^\prime$), and 
each of the faces have already been triangulated. In other words, we see that we have a topological $\mD^{k+1}$, along with
a given triangulation of $\partial \mD^{k+1}$. But it is now easy to extend: just cone the given triangulation on the
boundary inwards. Performing this process on each of the $\sigma ^k\times I$ now provides us with a triangulation
of the set $\mT_i^{k}\times I$. This results in a triangulation of the 3-sphere with the following two properties:
\begin{itemize}
\item the triangulation contains a collection of squares, whose union realize the given link $k$,
\item for each of the squares, the union of the simplices incident the the square form a regular neighborhood $\mD^2 \times S^1$,  
triangulated as in Step 1, and
\item all of these regular neighborhoods are pairwise disjoint.
\end{itemize}

\vskip 5pt

\noindent {\bf Step 3: Getting rid of all other squares.}

\vskip 5pt

At this stage, we have constructed a triangulation of $S^3$, which contains a collection of squares realizing
the given link $k$. However, there are still two problematic issues: our triangulation might not be flag, and
it might fail the isolated squares condition. The third step is to modify the triangulation in order to ensure these two
additional conditions. To fix some notation, we will keep using $N_i$ to denote the regular neighborhood of
the squares we are interested in keeping. Recall that each of these is topologically a solid torus 
$\mD^2\times S^1$, with triangulation combinatorially isomorphic to the triangulation given in Step 1. We
will first modify the given triangulation in the complement of the $N_i$, and subsequently change it within
the regions $N_i$. 

Let us denote by $X$ the closure of the complement of the union of the $N_i$. This is topologically a 
3-manifold with boundary, equipped with a triangulation (from the previous two
steps). Now the standard method of obtaining a flag triangulation is to take the barycentric subdivision
of a given triangulation. But unfortunately, this process creates lots of squares. Recently,
Przytycki and \'Swi{\polhk{a}}tkowski \cite{PS}, building on earlier work of Dranishnikov \cite{Dr}, 
have found a different subdivision process that takes a 3-dimensional 
simplicial complex and returns a subdivision of the complex that is flag {\it and has no squares}. 
For an arbitrary simplicial complex $Z$, we will denote by $Z^*$ the simplicial complex 
obtained by applying this procedure to $Z$. We modify the given triangulation
of $S^3$ in two stages: first we modify the triangulation in $X$, by replacing $X$ by $X^*$.
Next, we describe the extension of this triangulation into the various components $N_i$. For the original 
triangulation of each of the $N_i$, we see that the thirty-six tetrahedra are of one of two types: 
\begin{enumerate}
\item[(a)] twenty-four of them are the join of one of the interior vertices with a triangle on $\partial N_i$, and
\item[(b)] twelve of them are the join of one of the four edges on the core square with an edge on $\partial N_i$. 
\end{enumerate}
Now the subdivision $X^*$ restricts to a subdivision on each simplex in 
$\partial N_i$, which changes the simplicial complex $\partial N_i$ into $(\partial N_i)^*$. 
The effect of this subdivision on simplices in $\partial N_i$ is to subdivide each edge in $\partial N_i$ 
into two, and to replace each original triangle by the subdivision in Figure 2. We extend the 
subdivision $(\partial N_i)^*$ of $\partial N_i$ to a subdivision $N_i^\prime$ of the original $N_i$ in the most
natural way possible:

\begin{figure}
\label{graph}
\begin{center}
\includegraphics[width=2in, angle=0]{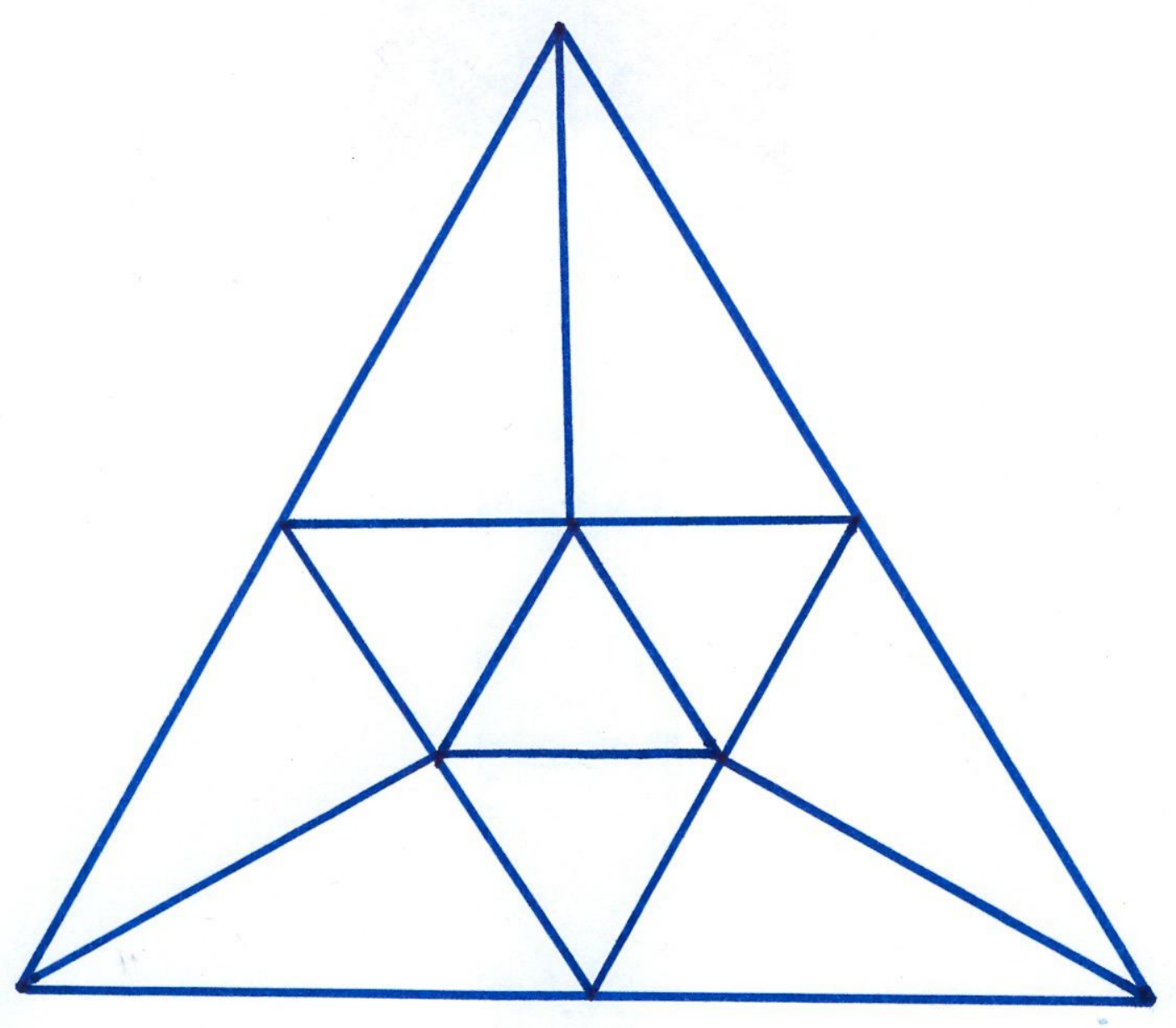}
\caption{Dranishnikov subdivision of triangles.}
\end{center}
\end{figure}

\begin{enumerate}
\item[(a)] each tetrahedron in $N_i$ that was a join of an interior vertex with a triangle $\sigma \subset \partial N_i$ gets
replaced by the join of the same vertex with $\sigma ^*$ (i.e. we cone the subdivision of $\sigma$ to the interior
vertex), subdividing the original tetrahedron into ten new tetrahedra (the cone over Figure 2), and
\item[(b)] each tetrahedron that was a join of an edge on the square with an edge on $\partial N_i$ gets replaced 
by two tetrahedra (i.e. the join of the internal edge with each of the two edges obtained from subdividing
the boundary edge).
\end{enumerate}
This changes the original triangulation on each $N_i$ into a new triangulation $N_i^\prime$ with a total of $264$ tetrahedrons.
We will continue to use the term {\it block} to refer to the subcomplexes of the $N_i^\prime$ that are subdivisions 
of the original blocks in $N_i$. Observe that, in each of the $N_i$, our subdivision process
did not introduce any new vertices in the interior of the $N_i$. As such, the core squares have been 
left unchanged (and we will still refer to them as the cores of the $N_i^\prime$).  

Finally, we note that by construction the two subdivisions $N_i^\prime$ of $N_i$, and $X^*$ of $X$ coincide on their common
subcomplex $\partial N_i = N_i\cap X$. In particular, they glue together to give a well defined triangulation 
$\Sigma$ of $S^3$.

\vskip 10pt

\noindent {\bf Step 4: Verifying that $\Sigma$ has the desired properties.}

\vskip 5pt

Note that the triangulation $\Sigma$ contains a copy of $X^*$, as well as copies 
of each $N_i^\prime$. These partition the triangulation $\Sigma$ into various pieces.

\begin{lemma} The complex $X^*$, the individual $N_i^\prime$, and the intersections $X^*\cap N_i^\prime$, are all 
full subcomplexes of $\Sigma$.
\end{lemma}

\begin{proof}
This follows easily from the following two facts:
\begin{itemize}
\item each of the intersections $X^* \cap N_i^\prime= (X \cap N_i) ^*$ is a full subcomplex of $X^*$, 
\item each of the intersections $X^*\cap N_i^\prime = \partial N_i^\prime$ is a full subcomplex of the corresponding $N_i^\prime$.
\end{itemize}
The first statement is a direct consequence of \cite[Lemma 2.10]{PS}, where it is shown that if $U$ is any subcomplex
of $W$, then $U^*$ is a full subcomplex of $W^*$. The second statement is a consequence of the construction
of the triangulation $N_i^\prime$, since by construction, each simplex of $N_i^\prime$ which is {\it not} contained in
$\partial N_i^\prime$ contains a vertex in the interior of $N_i^\prime$ (and hence in $N_i^\prime - \partial N_i^\prime$). 
\end{proof}

\begin{lemma}
The triangulation $N_i^\prime$ is flag.
\end{lemma}

\begin{proof}
Given a collection of pairwise incident vertices $V$, there are three possibilities: $V$ contains either
two, one, or no interior vertices of $N_i^\prime$. We consider each of these three cases in turn.

If $V$ contains no interior vertices, then $V\subset \partial N_i^\prime$, and since the latter is a full subcomplex of $N_i^\prime$
(see Lemma 4), $V$ is in fact a collection of vertices in $\partial N_i^\prime$ which are pairwise adjacent {\it within} 
$\partial N_i^\prime$. But recall that $\partial N_i^\prime$ is just the triangulation $(\partial N_i)^*$, 
hence is flag. This implies that $V$ spans out a simplex in $\partial N_i^\prime$.

If $V$ contains one interior vertex $v$, then, by the previous argument, $V-\{v\}$ spans a simplex in 
$\partial N_i^\prime = (\partial N_i)^*$ which is contained within some (maximal) 2-dimensional simplex $\sigma$
in $(\partial N_i)^*$. Note that, since all vertices $V-\{v\}$ are adjacent to the interior vertex $v$, they must lie
in the block $B$ corresponding to $v$. So the 2-dimensional simplex $\sigma \subset (\partial N_i)^*$ can 
additionally be chosen to lie within that same block $B$. This means that there exists a 2-dimensional simplex 
$\tau \in \partial N_i$ with the property that $\sigma$ is one of the 10 triangles in $\tau^*$ (see Figure 2). Finally,
observe $\tau$ must lie within the block $B$, so the join of $\tau$ with the interior vertex $v$ defines a tetrahedron
inside the original triangulation $N_i$ (of type (a) in the terminology of Step 3). But recall how the subdivision $(\partial N_i)^*$ 
of the triangulation $\partial N_i$ was extended into $N_i$: for tetrahedra of type (a), the subdivision on the
boundary was coned off to the interior vertex. This implies that the join of $\sigma$ and the vertex $v$ defines a 
tetrahedron in $N_i^\prime$, and as the set $V$ is a subset of the vertex set of this tetrahedron, we deduce that $V$ 
spans a simplex in $N_i^\prime$.

Finally, if $V$ contains two interior vertices $v,w$, let $B_v, B_w$ denote the corresponding blocks. Since $V-\{v,w\}$ is a 
collection of vertices in $\partial N_i^\prime= (\partial N_i)^*$ which are adjacent to {\it both} interior vertices, we see that the set $V-\{v,w\}$ 
must lie within $B_v\cap B_w$, which is a 1-dimensional complex homeomorphic to $S^1$ (subdivided into 6 consecutive 
edges). Since $V-\{v,w\}$ are pairwise adjacent, there is an edge $\sigma$ in $B_v\cap B_w$ whose vertex
set contains $V-\{v,w\}$. This edge is contained in a subdivision of an edge $\tau$ from the original triangulation 
$\partial N_i$, where $\tau$ is an edge which is common to the two blocks $B_v$ and $B_w$. In particular, the join $\omega*\tau$ of $\tau$
with the edge $\omega$ in the core joining $v$ to $w$ defines a tetrahedron in the original triangulation $N_i$ (of type (b) in the terminology
of Step 3). Again, from the way the subdivision $(\partial N_i)^*$ was extended inwards, we recall that the tetrahedra
$\omega * \tau$, being of type (b), gets replaced by two tetrahedra $\omega * \sigma$ and $\omega * \sigma^\prime$, where
$\tau^* = \sigma \cup \sigma^\prime$. Since the join of $\sigma$ and $\omega$ defines a tetrahedron in $N_i^\prime$, and the set
$V$ is a subset of the vertex set of this tetrahedron, we again deduce that $V$ spans a simplex in $N_i^\prime$.
\end{proof}

\begin{Cor} The triangulation $\Sigma$ is flag.
\end{Cor}

\begin{proof}
If all of the vertices are contained in $X^*$, then the claim 
follows immediately from the fact that $X^*$ itself is flag (see \cite[Proposition 2.13]{PS}). 
So we can now assume that at least one of the vertices is contained in the interior of one of the $N_i^\prime$. 

Note that an interior vertex in one of the $N_i^\prime$ has its closed star entirely contained within the same 
$N_i^\prime$. So we see that the tuple of pairwise adjacent vertices must be entirely contained within the 
same subcomplex $N_i^\prime$. But by Lemma 5, we have that each of the subdivided $N_i^\prime$ are themselves flag, 
finishing the proof.
\end{proof}

\begin{Prop} The only squares in $\Sigma$ are the cores of the various $N_i^\prime$.
\end{Prop}

\begin{proof}
To see this, let us start with an arbitrary square $(v_1, v_2, v_3, v_4)$ inside the triangulation
$\Sigma$. Our goal is to show that all four vertices must be interior vertices to a single $N_i^\prime$, which
would then force the square to be the core of the corresponding $N_i^\prime$. To this end, we first note that,
if the square does {\it not} contain any interior vertex to any of the $N_i^\prime$, then it is contained entirely
within $X^*$. But from Lemma 4, the latter is a full subcomplex of $\Sigma$, and by the result 
of Przytycki and \'Swi{\polhk{a}}tkowski \cite[Proposition 2.13]{PS}, has no squares. So we may assume that 
at least one of the vertices is an interior vertex to some $N_i^\prime$.  

If all the vertices are interior to $N_i^\prime$, then we are done, so by way of contradiction we can also assume 
that the square contains a vertex which is {\it not} interior to $N_i^\prime$ (which we will call {\it exterior} vertices 
to $N_i^\prime$). Now the square $(v_1, v_2, v_3, v_4)$ contains exactly four edges, and since it contains vertices which are both interior and 
exterior to $N_i^\prime$, we must have that at least two of the four edges must connect an interior vertex to an
exterior vertex (call these {\it intermediate} edges). 

We now argue that in fact the square must contain {\it exactly}
two intermediate edges. Indeed, if there were $\geq 3$ intermediate edges, then one could find a pair
of adjacent intermediate edges, which share a common exterior vertex. Up to cyclic relabeling, we may assume that $v_1$
is the exterior vertex. Considering the other endpoints
of these two intermediate edges, we see that $v_2, v_4$ are interior vertices
for $N_i^\prime$, which
are both adjacent to the exterior vertex $v_1 \in \partial N_i^\prime$. But this implies that the two blocks
whose bottoms contain $v_2$ and $v_4$ cannot be opposite, so must in fact be adjacent.
This forces $v_2$ and $v_4$ to be 
adjacent vertices in the core of $N_i^\prime$, contradicting the fact that $(v_1, v_2, v_3, v_4)$ forms a square.  
So our hypothetical square $(v_1, v_2, v_3, v_4)$ must have exactly two intermediate edges, leaving us 
with exactly two possibilities:
\begin{enumerate}
\item the intermediate edges are not adjacent in the square $(v_1, v_2, v_3, v_4)$, 
\item the intermediate edges are adjacent at an interior vertex of $N_i^\prime$, and the remaining edges are exterior.
\end{enumerate}
We now explain why each of these possibilities give rise to a contradiction.

\begin{figure}
\label{graph}
\begin{center}
\includegraphics[width=5in, angle=0]{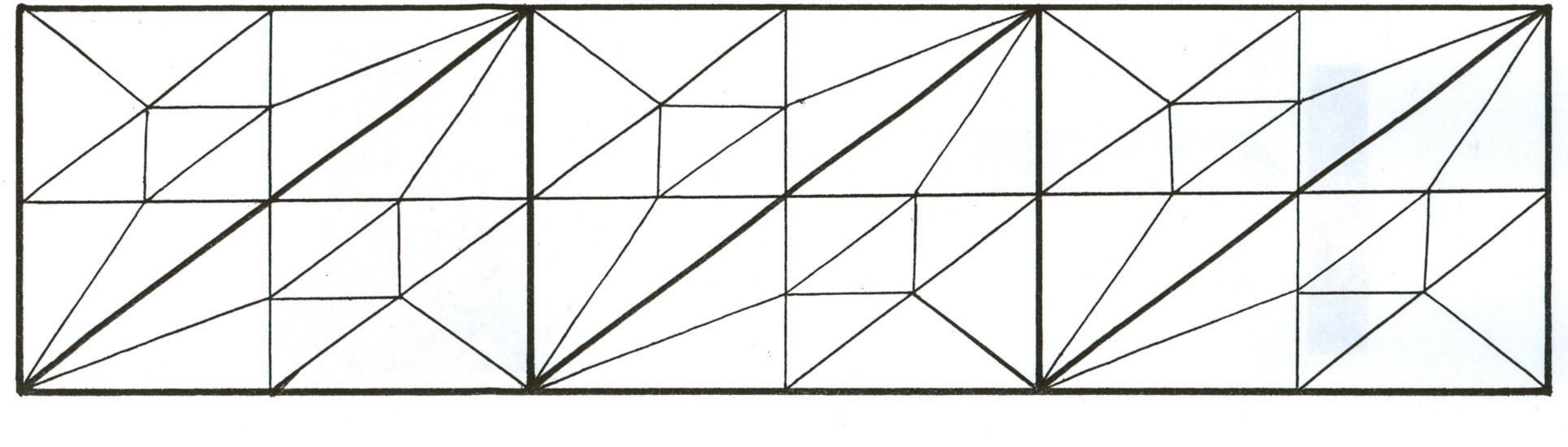}
\caption{Triangulation on the boundary of a block.}
\end{center}
\end{figure}

In case (1), we note that up to cyclic relabeling, we have that $v_1, v_2$ are adjacent vertices in the core of the $N_i^\prime$,
while $v_3, v_4$ are adjacent vertices in $\partial N_i^\prime$. We can also assume that the top of the block $B_1$ corresponding to
$v_1$ attaches to the bottom of the block $B_2$ corresponding to $v_2$. Now recall that an interior vertex is {\it only adjacent
to boundary vertices in its corresponding block}. Since $v_3$ is adjacent to $v_2$, we have that $v_3$ must lie in the 
block $B_2$. Similarly, the vertex $v_4$ being adjacent to $v_1$ must lie in the block $B_1$. Since $v_3$ and
$v_4$ are adjacent, we conclude that one of these two vertices must lie in the common boundary $B_1\cap B_2$. But
such a vertex is incident to both $v_1$ and $v_2$, violating the square condition for $(v_1, v_2, v_3, v_4)$.

It remains to rule out case (2). To this end, we may again assume that $v_1$ is the common interior vertex for the
two intermediate edges. Now if $B$ denotes the block corresponding to $v_1$, then we have that the boundary vertices
$v_2, v_4$, both being adjacent to $v_1$, must actually lie in $B$. Moreover, for $(v_1, v_2, v_3, v_4)$ to be a square,
we must have that $v_3$ is {\bf not} adjacent to $v_1$, and hence $v_3\notin B$. Since $v_3$ is adjacent to both the vertices
$v_2, v_4 \in B$, we see that the latter are either both in the top of $B$ or both in the bottom of $B$, while $v_3$ lies in
an adjacent block $B^\prime$. Let us assume that the vertices lie in the top of $B$ (the other case being completely analogous),
so that we can view $v_2,v_4$ as lying in the {\it bottom} of the block $B^\prime$. 

We now have the following situation occurring inside the boundary of the block $B^\prime$: we have two vertices $v_2, v_4$ lying 
in the bottom of the block, and we have a vertex $v_3$ which does {\bf not} lie in the bottom of $B^\prime$, but which is adjacent to both 
$v_2$ and $v_4$. Now recall that the triangulation of the block $B^\prime$ is a subdivision (given in Step 3) of a canonical
triangulation of the triangular prism. This subdivision takes the boundary of the original triangulation and applies the Dranishnikov
subdivision procedure to it: each edge gets subdivided into two, and each triangle gets replaced by the subdivision in Figure 2. 
The resulting triangulation on $S^1\times [0,1]$ is shown in Figure 3. In the illustration, the left and right side of the rectangle have to 
be identified,
and the ``bottom'' and ``top'' of the boundary of the block is precisely the bottom and the top of the rectangle. Note that this triangulation
actually consists of six original triangles (see Step 1), each of which has been subdivided into 10 triangles as in Figure 2 (see Step 3).
Finally, inspecting the
triangulation in Figure 3, we observe that there are exactly six vertices which are adjacent to two distinct 
vertices in the bottom of the block: these are the only possibilities for $v_3$. But for each of these six vertices, we see that the 
two adjacent vertices in the bottom of the block (i.e. the corresponding $v_2$ and $v_4$) are adjacent to each other, contradicting
the fact that $(v_1, v_2,v_3,v_4)$ was a square.

Since we've ruled out all other possibilities, we see that the square cannot contain {\it any} intermediate edges, i.e. the four vertices of our 
hypothetical square $(v_1,v_2,v_3,v_4)$ must all lie in the
interior of a single $N_i^\prime$. This implies that our square must coincide with the core of one of the $N_i^\prime$, as desired.

\end{proof}

It follows from Corollary 6 that the triangulation $\Sigma$ is flag, and from Proposition 7 that it has isolated squares with
type given by the original link $k$. This completes the proof of Theorem 3.

\section{Constructing the manifold.} 

In this section, we establish the Main Theorem. Our goal is to use some of the
triangulations of $S^3$ constructed in the previous section to produce a 4-dimensional 
manifold $M$ with the desired properties. In order to do this, we start by reviewing some
properties of the Davis complex for right angled Coxeter groups.

\vskip 10pt

Recall that one can associate to the 1-skeleton of {\it any} simplicial complex
$L$ a corresponding {\it right angled Coxeter group} $\Gamma_L$. This group has one generator
$x_i$ of order two for each vertex $v_i$ of the simplicial complex $L$, and a relation 
$x_ix_j=x_jx_i$ whenever the corresponding vertices $v_i, v_j$ are adjacent in $L$.
Let us consider the associated Davis complex $\tilde P_{L}$. This complex is obtained via 
the following procedure: we first consider the cubical complex $[-1,1]^{V(L)}$, that is to say, the 
standard cube with dimension equaling the number of vertices in the simplicial complex $L$. Now
every face of the cube is an affine translation of $[-1,1]^S$ for some subset $S\subset V(L)$, which
we call the {\it type} of the face. Consider the cubical subcomplex $P_L \subset [-1,1]^{V(L)}$ 
consisting of all faces whose type defines a simplex in $L$, and let $\tilde P_L$ to be its 
universal cover. Observe that the Coxeter group $\Gamma_L$ acts on $P_L$, where each generator
$x_i$ acts by reflection on the corresponding coordinate. The kernel of the resulting morphism $\Gamma_L
\rightarrow (\mZ _2)^{|V(L)|}$ coincides with the fundamental group of $P_L$. There is a natural
piecewise flat metric on $P_L$, obtained by making each $k$-dimensional face in the cubulation 
of $P_L$ isometric to $[-1,1]^k \subset \mR ^k$. 
Properties of the cubical complex $P_L$ are intimately related to properties of the simplicial complex
$L$. For instance, we have:

\begin{enumerate}
\item[(a)] if $L$ is a flag complex, then the piecewise flat metric on $P_L$ is locally CAT(0),
\item[(b)] the links of vertices in $P_L$ are canonically simplicially isomorphic to $L$,
\item[(c)] if $L$ is the join of two subcomplexes $L_1, L_2$, then the space $P_L$
splits isometrically as a product of $P_{L _1}$ and $P_{L _2}$,
\item[(d)] if $L ^\prime$ is a full subcomplex of $L$, then the natural inclusion induces a totally
geodesic embedding $P_{L^\prime}\hookrightarrow P_{L}$,
\item[(e)] if the geometric realization of $L$ is homeomorphic to an $(n-1)$-dimensional sphere, then 
$P_L$ is an $n$-dimensional manifold,
\item[(f)] if $L$ is a PL-triangulation of $S^{n-1}$ then $P_L$ is a PL-manifold, and
$\partial ^\infty \tilde P_L$ is homeomorphic to $S^{n-1}$,
\item[(g)] if $L$ is a {\it smooth} triangulation of $S^{n-1}$, then $P_L$ is a smooth manifold.
\end{enumerate}
These results are discussed in detail in the book \cite{Da1}. 

\vskip 10pt

In the previous section, we showed that given a prescribed link $k$ in $S^3$, one can construct a 
triangulation of $S^3$ with isolated squares, and with type the given link. Let us apply this result in the
special case where $k$ is a nontrivial knot inside $S^3$.  Let $L$ denote the corresponding
triangulation of $S^3$. Since we are in the special case of dimension $=3$, the triangulation 
$L$, in addition to being flag, is automatically PL and smooth. We now consider the cubical complex 
$M:= P_L$ associated to the corresponding right angled Coxeter group $\Gamma_L$. In view of our 
earlier discussion, we have the following:

\vskip 5pt

\noindent {\bf Fact 1:} The space $M$ is a smooth 4-manifold (from (g) above), and the natural
piecewise Euclidean metric on $M$ induced from the cubulation is locally CAT(0) (from (a) above). 
Furthermore, the boundary at infinity of $\tilde M$ is homeomorphic to $S^3$ (from (f) above), and 
$\tilde M$ is diffeomorphic to $\mR^4$.

\vskip 10pt

The very last statement in {\bf Fact 1} can be deduced from work of Stone (see \cite[Theorem 1]{St}),
who showed that a metric (piecewise flat) polyhedral complex which is both CAT(0) and a PL-manifold 
without boundary must in fact be PL-homeomorphic to the appropriate $\mathbb R^n$. Since our $\tilde M$ satisfies
these conditions, this ensures that $\tilde M$ is PL-homeomorphic to the standard $\mathbb R^4$.
But in the 4-dimensional setting, there is no difference between PL and smooth, so
$\tilde M$ is in fact diffeomorphic to $\mathbb R^4$.

Our goal is now to show that $M$ has the
properties postulated in our Main Theorem. Note that properties (1) and (2) are included in {\bf Fact 1},
while property (4) can be easily deduced from property (3) (see the comment after the proof of Proposition 1).
So we are left with establishing property (3): that $\pi_1(M)$ cannot be isomorphic
to the fundamental group of any nonpositively curved Riemannian manifold.
This last property will be established by looking at the large scale geometry of flats inside
the universal cover $\tilde M$. 

As a starting point, let us describe some flats inside $\tilde M$. Observe that each square inside
the triangulation $L$ is a full subcomplex isomorphic to a 4-cycle $\square$. The right angled
Coxeter group associated to a 4-cycle is a direct product of two infinite dihedral groups 
$\Gamma_\square \cong D_\infty \times D_\infty = (\mZ _2 * \mZ_2) \times (\mZ_2 * \mZ_2)$ 
(see (c) above).
The corresponding complex $P_\square$ is isometric to a flat torus (with cubulation given by 16
squares, obtained via the identification $S^1\times S^1 = \square  \times \square$).
By considering the unique square inside the triangulation $L$, we obtain:

\vskip 10pt

\noindent {\bf Fact 2:} $M $ contains a totally geodesic $2$-dimensional flat torus $T^2$ (see (d) above). 
Furthermore, at any vertex $v\in T^2 \subset M$ of the cubulation, we have that the torus $T^2$ is {\it locally knotted} inside
the ambient $4$-dimensional manifold $M$ (see (b) above), in that there is a canonical simplicial isomorphism
$\big(lk_v(M), lk_v(T^2)\big) \cong (L, k)$ where $k$ is the unique (knotted) square in the triangulation $L$.

\vskip 10pt

Since the embedding $T^2 \hookrightarrow M$ is totally geodesic, by lifting to the universal cover, we obtain
a $2$-dimensional flat $F \hookrightarrow \tilde M$ which is locally knotted at lifts of vertices. 
This induces an embedding of the corresponding boundaries at infinity, giving us an embedding of 
$\partial ^\infty F \cong S^1$ into $\partial ^\infty \tilde M \cong S^3$. 
The rest of our argument will rely on the following ``local-to-global'' assertion:

\vskip 10pt

\noindent {\bf Assertion:} The embedding $\partial ^\infty F \cong S^1$ into $\partial ^\infty \tilde 
M \cong S^3$ defines a nontrivial knot in the boundary at infinity of $\tilde M$.

\vskip 10pt

That is to say, the ``local knottedness'' of the flat propagates to ``global knottedness'' of its boundary at infinity.
For the sake of exposition, we delay the proof of the assertion, and first show how we can use it to deduce the 
Main Theorem. To this end, let us assume that $(M^\prime, g)$ is a closed manifold equipped with a Riemannian
metric of nonpositive sectional curvature, and that we are given an isomorphism of fundamental groups 
$\phi:\Gamma= \pi_1(M) \rightarrow \pi_1(M^\prime)$. From this assumption, we want to work towards a contradiction.

\vskip 10pt

The first step is to use the isomorphism 
of fundamental groups to obtain an equivariant homeomorphism between the corresponding boundaries at infinity.
As a cautionary remark, we recall that given a pair $X_1,X_2$ of CAT(0)-spaces with 
geometric $G$-actions, a celebrated example of Croke and Kleiner \cite{CK} shows that the corresponding 
boundaries at infinity $\partial ^\infty X_1$ and $\partial ^\infty X_2$ need {\bf not} be homeomorphic. Even if 
the boundaries at infinity {\it are} homeomorphic, an example of Buyalo \cite{Bu} shows that the homeomorphism might
{\bf not} be equivariant with respect to the $G$-action. 

In his thesis \cite{H}, Hruska introduced CAT(0)-spaces with {\it isolated flats}. Subsequent work of Hruska and 
Kleiner \cite{HK} established the following two foundational results for CAT(0)-spaces with isolated flats:
\begin{enumerate}
\item for a pair $X_1, X_2$ of CAT(0)-spaces with geometric $G$-actions, if
$X_1$ has isolated flats, then so does $X_2$ (see \cite[Corollary 4.1.3]{HK}), and there is a $G$-equivariant
homeomorphism between $\partial ^\infty X_1$ and $\partial ^\infty X_2$ (see \cite[Theorem 4.1.8]{HK}).
\item for a group $G$ acting geometrically on a CAT(0)-space $X$, we have that $X$ has the isolated flats property
if and only if $G$ is a relatively hyperbolic group with respect to a collection of virtually abelian subgroups of rank
$\geq 2$ (see \cite[Theorem 1.2.1]{HK}).
\end{enumerate}
As such, if we could establish that our {\it group} $\Gamma$ is a relatively hyperbolic group with respect to a collection of
virtually abelian subgroups of rank $\geq 2$, then result (2) above would ensure that our CAT(0)-manifold $\tilde M$
has the isolated flats property. Result (1) above would then give the desired $\Gamma$-equivariant homeomorphism 
between $\partial ^\infty \tilde M$ and $\partial ^\infty \tilde M^\prime$. So our next goal is to establish:

\vskip 10pt

\noindent {\bf Fact 3:} The group $\Gamma=\pi_1(M)$ is hyperbolic relative to the collection of all virtually abelian subgroups
of $\Gamma$ of rank $\geq 2$.

\vskip 10pt

The notion of a group $G$ being relatively hyperbolic with respect to a collection $\mathcal A$ of subgroups of $G$ was
originally suggested by Gromov \cite{Gr}, whose approach was later formalized by Bowditch \cite{Bo}. Alternate formulations
appear in Farb's thesis \cite{Fa}, in work of Dru\c tu and Sapir \cite{DrSa}, and in the memoir of Osin \cite{Os}. 
We refer the reader to the original sources
for a detailed definition as well as basic properties of such groups. For our purposes, we merely need to know that the 
property of a group $G$ being hyperbolic relative to a collection of virtually abelian subgroups of rank $\geq 2$ 
is inherited by finite index subgroups of $G$. In particular, to show the desired property for $\Gamma$, we see that it is 
sufficient to establish that our
original Coxeter group $\Gamma_L$ is relatively hyperbolic with respect to higher rank virtually abelian subgroups 
(since $\Gamma \leq \Gamma_L$ is of finite index). 

Caprace \cite[Cor.~D~(ii)]{Ca} recently provided a criterion for deciding whether a Coxeter group is 
hyperbolic relative to the collection of its higher rank virtually abelian subgroups. In the right-angled case the condition is that the flag complex $L$ which defines $\Gamma_L$ contains no full subcomplex isomorphic to the suspension $\Sigma K$ of a subcomplex $K$ with 3 vertices which is either
\begin{enumerate}
\item[(a)]
the disjoint union of 3 points, or
\item[(b)]
the disjoint union of an edge and 1 point.
\end{enumerate}
In both cases $\Sigma K$ does not have isolated squares.  Since the Coxeter group $\Gamma_L$ with which we are working is associated to a triangulation $L$ of $S^3$ with isolated squares, we conclude that $\Gamma_L$ is relatively hyperbolic with respect to  the collection of all virtually abelian subgroups of rank $\geq 2$.  Hence, {\bf Fact 3}.

\vskip 10pt 

Applying Hruska and Kleiner's results from \cite{HK}, we conclude that the original $\tilde M$ is a CAT(0)-space
with the isolated flats property, and that there exists a $\Gamma$-equivariant homeomorphism from $\partial ^\infty \tilde
M$ to $\partial ^\infty \tilde M^\prime$. The nontrivial knot $\partial ^\infty F \cong S^1$ inside $\partial
^\infty \tilde M \cong S^3$ appearing in the {\bf Assertion} can be identified with the limit set of the corresponding
subgroup $\pi_1(T^2) \cong \mZ^2 \leq \Gamma = \pi_1(M)$. Since we have an equivariant homeomorphism between 
the boundaries at infinity of $\tilde M$ and $\tilde M^\prime$, this immediately yields:

\vskip 10pt

\noindent {\bf Fact 4:} The boundary at infinity $\partial ^\infty \tilde M^\prime$ is homeomorphic to $S^3$, and the 
limit set of the canonical $\mZ^2$-subgroup in $\Gamma\cong \pi_1(M^\prime)$ defines a nontrivial knot  
$S^1\hookrightarrow \partial ^\infty \tilde M^\prime\cong S^3$.

\vskip 10pt

On the other hand, the flat torus theorem implies that there exists a $\mZ^2$-periodic flat 
$F^\prime \hookrightarrow \tilde M^\prime$, with the property that $\partial ^\infty F^\prime$ coincides with
the limit set of the $\mZ^2$. In particular, $\partial ^\infty F^\prime$ defines a nontrivial knot inside $\partial ^\infty M^\prime$.
But taking any point $p\in F^\prime$, we note that geodesic retraction provides a homeomorphism 
$\rho: \partial ^\infty \tilde M^\prime \rightarrow T_p\tilde M^\prime$. This homeomorphism takes the knotted
subset $\partial ^\infty F^\prime$ lying inside $S^3 \cong \partial ^\infty M^\prime$ to the {\it unknotted}
subset $T_pF^\prime$ lying inside $S^3\cong T_p\tilde M^\prime$. This contradiction
allows us to conclude that no such Riemannian manifold $(M^\prime, g)$ can exist.

\vskip 20pt

So in order to complete the proof of the Main Theorem, we are left with establishing the {\bf Assertion}. We note that
a similar result was shown in the setting of CAT(-1)-manifolds by Farrell and Lafont \cite{FL}, the proof of which extends
almost verbatim to yield the {\bf Assertion}. For the convenience of the reader,
we provide a (slightly different) self-contained argument for the {\bf Assertion}. 

The basic idea is as
follows: picking a vertex $v \in F$, we have a geodesic retraction map $\rho:\partial ^\infty \tilde M\rightarrow lk_v(\tilde M)$.
Under this map, we see that $\partial ^\infty F$ maps to the link $lk_v (F)$ inside $lk_v (\tilde M)$. But recall from 
{\bf Fact 2} that the torus is locally knotted in $\tilde M$, i.e. the pair $\big( lk_v(\tilde M), lk_v(F) \big)$ is simplicially 
isomorphic to $(S^3, k)$, where $S^3$ is the 3-sphere equipped with the triangulation $L$, and $k$ is the 
knot in $S^3$ given by the unique square in the triangulation $L$. Now the retraction map $\rho$ is {\it not}
a homeomorphism, but is nevertheless ``close enough'' to a homeomorphism for us to use it to compare the pair
$\big(\partial ^\infty \tilde M, \partial ^\infty F \big)$ with the knotted pair $\big( lk_v(\tilde M), lk_v(F) \big) \cong (S^3, k)$.
More precisely, for any given subset $Z\subset lk_v(\tilde M)\cong S^3$ we denote by $Z_\infty$ the corresponding pre-image
$Z_\infty:= \rho^{-1}(Z)$ inside $\partial ^\infty \tilde M$. Then we have:

\vskip 10pt

\noindent{\bf Fact 5:} \cite[Proposition 2, pg. 627]{FL} For any open set $U \subset lk_v(\tilde M)$, the map $\rho:U_\infty\rightarrow U$ 
is a proper homotopy equivalence. Moreover, the map $\rho$ is a {\it near-homeomorphism}, i.e. can be approximated arbitrarily 
closely by homeomorphisms.

\vskip 10pt

This is shown by identifying $U_\infty$ with the inverse limit of the sets $\{U_r\}_{r\in \mR^+}$, where each $U_r$ is the pre-image of
$U$ under the geodesic projection from the sphere $S_v(r)$ of radius $r$ centered at $v$ to the link at $v$. For $r>s$, the bonding maps 
$\rho_{r,s}: U_r \rightarrow U_s$ are given by geodesic retraction, and the canonical map $\rho_{\infty, s}$ from $U_\infty = \varprojlim \{U_r\}$ 
to each individual $U_s$ coincides with the geodesic retraction map. Since the link $lk_v(\tilde M)$ can be identified with $S_\epsilon(r)$,
a small enough $\epsilon$-sphere centered at $v$, the map $\rho$ can be identified with the canonical map $\rho_{\infty, \epsilon}$ from 
$U_\infty = \varprojlim \{U_r\}$ to the corresponding $U_\epsilon = U$. Now by results of Davis and Januszkiewicz \cite[Section 3]{DJ} each
of the bonding maps $\rho_{r,s}$ are cell-like maps, i.e. point pre-images have the shape of a point (see Dydak and Segal \cite{DySe} for 
background on shape theory). Since the shape functor commutes with inverse limits, and since $\rho = \rho_{\infty, \epsilon}$, we see 
that $\rho$ is also a cell-like map. A result of Edwards  \cite[Section 4]{Ed} now implies that $\rho$ is a proper homotopy equivalence,
while work of Armentrout \cite{Ar} ensures that $\rho$ is a near homeomorphism.

\vskip 10pt

Now to show that $\partial ^\infty F$ defines a nontrivial knot in $\partial ^\infty \tilde M$, we need to establish that 
the complement $\partial ^\infty \tilde M - \partial ^\infty F$ cannot be homeomorphic to $S^1\times \mR^2$. This will follow
if we can show that $\pi_1\big(  \partial ^\infty \tilde M - \partial ^\infty F \big)$ is a non-abelian group. To do this, let us 
decompose $\partial ^\infty \tilde M - \partial ^\infty F$ into a union of a suitable pair of open sets. We start by decomposing
$lk_v(\tilde M)$, and will then use the map $\rho$ to ``lift'' this decomposition to $\partial ^\infty \tilde M$.
Let  $lk_v(F) \subset N_1\subset N_2 \subset lk_v(\tilde M)$ be nested open regular neighborhoods of the knot $k=lk_v(F)$ inside 
$S^3 \cong lk_v(\tilde M)$. Define open sets in $lk_v(\tilde M)$ by setting $U_2:= N_2$, and $U_1:= lk_v(\tilde M) - \bar N_1$, 
where $\bar N_1$ denotes the closure of $N_1$. Note that we have homeomorphisms $U_2\cong S_1\times \mD ^2$ and 
$U_1\cap U_2 \cong N_2 - \bar N_1\cong S^1\times S^1\times \mR$, while $U_1$ is homeomorphic to the complement of the
nontrivial knot $k\subset S^3$. So at the level of $\pi_1$, we have that (a) $\pi_1(U_1\cap U_2) \cong \mZ \oplus \mZ$, and
(b) $\pi_1(U_1)$ is a non-abelian group. The latter fact follows from work of Papakyriakopoulos \cite{Pa}, who
showed that $\pi_1$ of the complement of a nontrivial knot cannot be isomorphic to $\mZ$. But by Alexander duality such a group 
must have abelianization isomorphic to $\mZ$, hence cannot be abelian. 

Now corresponding to this decomposition of $lk_v(\tilde M)$, we have an associated open decomposition of $\partial ^\infty \tilde M$ in
terms of the corresponding $(U_1)_\infty$, $(U_2)_\infty$. We now define an open decomposition of 
$\partial ^\infty \tilde M - \partial ^\infty F$ by setting $U:= (U_1)_\infty$ and $V:= (U_2)_\infty - \partial ^\infty F$. The 
intersection satisfies $U\cap V = (U_1\cap U_2) _\infty$. Applying {\bf Fact 5} to the discussion in the previous paragraph,
we obtain that (a) $\pi_1(U\cap V) \cong \mZ\oplus \mZ$, and (b) $\pi_1(U)$ is non-abelian. From Seifert-Van Kampen, we have:
$$ \pi_1\big(\partial ^\infty \tilde M - \partial ^\infty F \big) = \pi_1(U) *_{\pi_1(U\cap V)} \pi_1(V)$$
So to see that $\pi_1\big(\partial ^\infty \tilde M - \partial ^\infty F \big)$ is non-abelian, it suffices to show that the non-abelian
group $\pi_1(U)$ injects into the amalgamation. But this will follow from:

\vskip 10pt

\noindent {\bf Fact 6:} The map $i_*:\pi_1(U\cap V) \rightarrow \pi_1(V)$ induced by inclusion is injective. 

\vskip 10pt

To establish {\bf Fact 6}, we first choose a suitable basis for $\pi_1(U\cap V)\cong \mZ \oplus \mZ$. Recall that the map $\rho$
gives a proper homotopy equivalence between $U\cap V = (U_1\cap U_2)_\infty$ and the space $U_1\cap U_2 = N_2 - \bar N_1$,
where $N_1\subset N_2$ are nested open regular neighborhoods of the knot $k$. Since $U_2 = N_2$ can be identified with 
$S^1\times \mD^2$, where $S^1\times \{0\}$ corresponds to the knot $k$, we choose the generators for 
$\pi_1(N_2 - \bar N_1)\cong \mZ \oplus \mZ$ to have the following two properties:
\begin{enumerate}
\item[\bf{(A)}] the generator $\langle 1, 0\rangle$ maps to a generator represented by $[S^1\times \{0\}] \in \pi_1(N_2)\cong \mZ$ 
under the obvious inclusion, and 
\item[\bf{(B)}] the generator $\langle 0, 1\rangle$ is chosen so that a representative curve exists which, under the natural inclusion 
into $N_2\cong S^1\times \mD^2$, projects to a generator for $\pi_1(\mD^2 -\{0\}) \cong \mZ$ in the $\mD^2$-factor, and is null-homotopic
in $N_2$.
\end{enumerate}
We choose the generators of $\pi_1(U\cap V)\cong \mZ \oplus \mZ$ to map to the above two generators of $\pi_1(U_1\cap U_2)$ 
under the homotopy equivalence $\rho$. 

To verify that $i_*: \pi_1(U\cap V) \rightarrow \pi_1(V)$ is injective, we first argue that an element $\langle a, b\rangle \in 
\ker (i_*)$ must satisfy $a=0$. Consider the commutative diagram:
$$\xymatrix{
\mZ \oplus \mZ \cong \pi_1(U\cap V) \ar[d]^{\rho_*}  \ar[r]^-{i_*} & \pi_1(V) \ar[r] & \pi_1\big( (N_2)_\infty\big) \cong \mZ \ar[d]^{\rho_*} \\
\mZ \oplus \mZ \cong \pi_1(U_1\cap U_2) \ar[rr] &  &  \pi_1( N_2) \cong \mZ \\
}
$$
where all horizontal arrows are induced by the obvious inclusions, and the two vertical arrows are the isomorphisms induced
by the geodesic retraction maps. By the choice of the basis on $\pi_1(U\cap V)$, we have that $\rho_*(\langle a, b\rangle)=
\langle a, b\rangle \in \pi_1(U_1\cap U_2)$, which by property {\bf (A)} maps to $a\in \mZ \cong \pi_1(N_2)$. From the commutativity
of the diagram, we conclude that if $\langle a, b\rangle \in \ker (i_*)$, then $a=0$. Our next goal is to show that $b=0$.

Given a pair $\eta_1, \eta_2$ of disjoint oriented curves in $S^1\times \mD^2$, 
with $\eta_1$ null-homotopic, there is a well-defined linking number $L(\eta_1,\eta_2)$. For smooth curves this is obtained by
looking at the oriented intersection number of $\eta_2$ with a smooth bounding disk for the curve $\eta_1$, and for continuous 
curves one uses an approximation by smooth curves. This linking number has the property that if $\eta_1\sim \eta_1^\prime$ 
(respectively $\eta_2\sim \eta_2^\prime$) are two curves
homotopic to each other {\it in the complement of $\eta_2$} (respectively $\eta_1$), then 
$L(\eta_1, \eta_2^\prime)=L(\eta_1, \eta_2)=L(\eta_1^\prime, \eta_2)$. 

Now from the choice of basis on $\pi_1(U\cap V)$, 
along with property {\bf (B)}, we can choose a representative curve $\gamma$ for the element $\langle 0,b\rangle \in \ker (i_*) 
\subset \pi_1(U\cap V)$ with the property that the image curve $\rho(\gamma)
\subset U_1\cap U_2 \subset N_2 \cong S^1\times \mD^2$ projects to $b$ times a generator for $\pi_1(\mD^2- \{0\})$. 
One can easily check that this forces $L\big(\rho(\gamma) , S^1\times \{0\}\big) = \pm b$. Applying {\bf Fact 5}, we can find 
a homeomorphism $\rho^\prime : (N_2)_\infty \rightarrow N_2$ which is $\epsilon$-close to the map $\rho$. In view of the discussion
above, and recalling that the curve $S^1\times \{0\}$ corresponds to $lk_v(F)= \rho(\partial ^\infty F)$, this gives us that:
\begin{eqnarray*}
\pm b &=& L\big(\rho(\gamma) , S^1\times \{0\}\big) = L\big(\rho(\gamma), \rho(\partial ^\infty F) \big)\\
&=& L\big(\rho^\prime (\gamma), \rho(\partial ^\infty F)\big) = L\big(\rho^\prime (\gamma), \rho^\prime(\partial ^\infty F) \big) = L\big(\gamma, \partial ^\infty F \big) \\
\end{eqnarray*}
\vskip -15pt
\noindent where for the last equality, we use the fact that $\rho^\prime$ is a homeomorphism, and hence preserves the linking 
number. But since $\langle 0, b\rangle \in \ker (i_*)$, we also have that $\gamma$ bounds a disk in $V = (N_2)_\infty- \partial ^\infty F$, 
which implies that $L\big(\gamma, \partial ^\infty F \big)=0$. This now forces $b=0$, completing the proof of
{\bf Fact 6}. 

\vskip 10pt

Since the non-abelian group $\pi_1(U)$ injects into 
$\pi_1\big(\partial ^\infty \tilde M - \partial ^\infty F \big)$, we obtain that 
$\partial^\infty F \cong S^1$ defines a nontrivial knot in $\partial ^\infty \tilde M \cong S^3$, establishing
the {\bf Assertion}, and finishing off the proof of the Main Theorem.

\section{Concluding remarks.}

Finally, we point out a few interesting questions that come up naturally from this work. As discussed in Section 2.2, locally 
CAT(0)-manifolds whose universal covers are {\bf not} diffeomorphic to $\mR^n$ cannot support a Riemannian 
smoothing. In dimensions $n\neq 4$, there is no difference between ``homeomorphic to $\mR^n$'' and
``diffeomorphic to $\mR^n$''. In contrast, it is known that $\mR^4$ supports many distinct smooth structures (in
fact, continuum many). Moreover, the method used to construct the Davis examples of closed aspherical manifolds
whose universal covers are not homeomorphic to $\mR^n$ requires $n\geq 5$.  So one can ask:

\vskip 10pt

\noindent {\bf Question:} Can one find locally CAT(0) closed 4-manifolds $M^4$ with the property that
their universal covers $\tilde M^4$ are
\begin{enumerate}
\item not homeomorphic to $\mR^4$?
\item homeomorphic, but not diffeomorphic to $\mR^4$?
\end{enumerate}

\vskip 10pt

\noindent 
Paul Thurston \cite{Th} proved that $\tilde M^4$ must be homeomorphic to $\mR^4$ if it has at least one ``tame'' point.  
We remark that the result of Stone \cite{St} tells us that there is no hope of constructing 
such examples via piecewise flat metric complexes (for their universal covers would then have
to be diffeomorphic to the standard $\mathbb R^4$). Moreover, if one asks instead for {\it aspherical} 
closed $4$-manifolds, we remark that Davis \cite{Da2} has constructed examples where the universal
cover is {\it not} homeomorphic to $\mathbb R^4$ (but it is unknown whether those examples support
a locally CAT(0)-metric). 

\vskip 5pt

Now concerning the dimension restriction in our construction, we note that this was due to the need for 
finding triangulations of spheres with the property that the associated Davis complex had the isolated
flats condition (in order to obtain a well-defined boundary at infinity). The ``isolated squares'' condition we 
introduced was designed to ensure that Caprace's criterion was fulfilled. Attempting to generalize this 
construction to higher dimensions, the difficulty we run into is that, by work of Januszkiewicz and \'Swi{\polhk{a}}tkowski
\cite[Section 2.2]{JS} (see also the discussion in \cite[Appendix]{PS}), there is no higher-dimensional analogue 
of the Dranishnikov-Przytycki-\'Swi{\polhk{a}}tkowski procedure for 
modifying triangulations in order to get rid of squares. 

%On the other hand, Caprace's criterion is quite 
%flexible, in that it is formulated in terms of suspensions (rather than the special case of squares). This 
%motivates the:
%
%\vskip 10pt
%
%\noindent {\bf Question:} Can one find flag triangulations of $S^n$ ($n\geq 4$) so that the associated
%Davis complex is CAT(0) with the isolated flats property? Can one find such examples where the Davis
%complex contains a codimension two flat which is knotted at infinity?
%
%\vskip 10pt
%
%\todo{Previous question needs to be tweaked.}

Finally, we remark that our construction relies on the presence of flats with specific
large scale behavior in order to obstruct Riemannian smoothings. As such, our methods require 
the presence of zero curvature. If one desires examples which are {\it strictly negatively curved}, 
we are brought to the following:

\vskip 10pt

\noindent {\bf Question:} Can one construct examples of smooth, locally CAT(-1)-manifolds $M^n$ with the property 
that $\partial ^\infty \tilde M$ is homeomorphic to $S^{n-1}$, but which do {\it not} support any Riemannian
metric of nonpositive sectional curvature?

\end{document}